%
\magnification=\magstep1
\input amstex
\UseAMSsymbols
\input pictex
\NoBlackBoxes
  \font\gross=cmbx10 scaled\magstep1 
   \font\rmk=cmr8    \font\itk=cmti8  \font\ttk=cmtt8

   \newcount\notenumber
   
   \def\note{\advance\notenumber by 1 
       \plainfootnote{$^{\the\notenumber}$}}

\def\Ext{\operatorname{Ext}}

\def\Rahmen#1%
   {$$\vbox{\hrule\hbox%
                  {\vrule%
                       \hskip0.5cm%
                            \vbox{\vskip0.3cm\relax%
                               \hbox{$\displaystyle{#1}$}%
                                  \vskip0.3cm}%
                       \hskip0.5cm%
                  \vrule}%
           \hrule}$$}

   
\vglue2truecm

\centerline{\gross Exceptional representations of quivers.}     
		   		       \bigskip
\centerline{Claus Michael Ringel}
		  	    \bigskip\medskip 
\plainfootnote{}
{\rmk 2010 \itk Mathematics Subject Classification. \rmk 
Primary 
        16G20, 
        16G60. 
Secondary:
        16D90, 
       	17B22.   
}
{\narrower\narrower Abstract. 
Let $Q$ be a connected directed quiver with $n$ vertices. 
We show that $Q$ is representation-infinite if and only if there do exist $n$ isomorphism classes
of exceptional modules of some fixed length $t\ge 2.$
\par}
	\bigskip\bigskip 

Let $Q$ be a connected directed quiver with $n$ vertices. 
A (finite-dimensional) representation $M$ of $Q$ is said to be {\it exceptional}
provided $M$ is indecomposable and $\Ext^1(M,M) = 0.$ Let $e(t) = e_Q(t)$ be the number of isomorphism
classes of exceptional representations of length $t$. Of course, we have $e(1) = n$ and if
$Q$ has no multiple arrows, then $e(2)$ is the number of arrows of $Q$.
	\medskip
{\bf Theorem.} {\it Let $Q$ be a connected directed quiver with $n$ vertices and assume that $Q$
has no multiple arrows. Then the following conditions are equivalent:}
\item{(i)} {\it The quiver $Q$ is representation-infinite.}
\item{(ii)} {\it There do exist $n$ isomorphism classes of 
    indecomposable representations of some fixed length $t\ge 2$.}
\item{(iii)} {\it There is $t$ with $2 \le t \le \min(n,7)$ such that $e(t) \ge n$, 
    $e(s) = n-1$ for $2\le s < t$ and such that 
    all indecomposable representations of length at most $t$
    are exceptional.}
	\medskip
\noindent
Note that according to (iii), there is $t\ge 2$ with $e(t)\ge n$. The smallest such 
number $t$ is an invariant of the quiver $Q$ and may be denoted by $t_Q = t.$
	\bigskip 
Proof of Theorem. The implication (ii) $\implies$ (i) is a direct consequence of an old result of 
Kostant (1959) which asserts that for a finite root system of rank $n$ and $t\ge 2$, there are 
at most $n-1$ positive roots of height $t$, see [K] or also [H], Theorem 3.20. Namely, if $Q$ is
representation-finite, then according to Gabriel [G] there is a bijection between the
isomorphism classes of the indecomposable representations of $Q$ of length $t$ and the 
positive roots (of a corresponding root system) of height $t$. The implication 
(iii) $\implies$ (ii) is trivial. Thus it remains to show the implication
(i) $\implies$ (iii).

Thus, we assume that $Q$ is representation-infinite. We will see that there is some natural number
$t\ge 2$ with $e_Q(t) \ge n$; the smallest such number $t$ will be denoted by $t_Q$.
If $Q$ is not a tree, then the number of arrows
is at least $n$, thus there are at least $n$ isomorphism classes of indecomposable modules
of length $2$. Since there are no multiple arrows, 
all indecomposable representations of length at most $2$ are exceptional. In particular, 
we see that $e(2) \ge n$ and therefore $t_Q = 2.$ Altogether we have shown that the assertion
(iii) with $t=2$ holds.
	\medskip 
We therefore can assume from now on 
that $Q$ is a tree, thus $e(2) = n-1.$ We insert the following lemma.
	\medskip
{\bf Lemma.} 
{\it Let $Q'$ be a subquiver of $Q$ with $n'$ vertices. If $2\le t\le n',$ then $e_Q(t) \ge e_{Q'}(t) + n-n'$.}
	\medskip
Proof of Lemma. It is sufficient to consider the case $n' = n-1$, thus $Q'$ is obtained from $Q$
by deleting one vertex, say $\omega$ (and all the arrows involving $\omega$). We consider the $e(t)$ 
exceptional
representations of $Q'$ as representations of $Q$; they are exceptional representations of
$Q$. Since $t \le n$, there is a connected subquiver $Q''$ 
of $Q$ which contains the vertex $\omega$ and such that $Q''$ has
precisely $t$ vertices. Since $Q$ is a tree, also  
$Q''$ is a tree, thus there is a unique thin indecomposable representation $M$ of $Q$ with
support equal to $Q''$ 
 (a representation is said to be 
{\it thin} provided all Jordan-H\"older multiplicities are bounded by $1$). 
Of course, a thin indecomposable representation of a tree quiver is exceptional. 
This shows that $e_Q(t) \ge e_{Q'}(t)+1 = 
e_{Q'}(t)+n-n'$. 
 
	\medskip 
We return to the proof of Theorem. 
Assume that $Q$ has a subquiver $Q'$ of type $\widetilde{\Bbb D}_m$ with $m\ge 4.$
Now $Q'$ has a subquiver $Q''$ of type $\Bbb A_{m-1}$, let $i,j$ be the two vertices 
of $Q'$ which do not belong to $Q''.$ Then $Q''$ has precisely $m-3$ connected subquivers with $3$
vertices. The quiver $Q'$ has additional connected subquivers with 3 vertices: for $m=4$
there are $5$ such subquivers, for $m\ge 5$ only $4$. Thus, we see that
$e_{Q'}(3) = (m-3)+5 = m+2$, whereas $e_{Q'}(3) = (m-3)+4 = m+1$, in case $m \ge 5$. 
Taking into account the Lemma, we see that $e_Q(3) \ge n+1$ in case $Q$ has a subquiver 
of type $\widetilde{\Bbb D}_4$, and 
that $e_Q(3) \ge n$ in case $Q$ has a subquiver 
of type $\widetilde{\Bbb D}_m$, with $m\ge 5.$ In both cases, we see that $t_Q = 3$. 
	
Next, we look at the quivers $Q$ with precisely one branching vertex and 
precisely $3$ arms. In this case, we claim that $e(3) = n-1.$ The proof is
by induction on $n$, starting with a quiver of type $\Bbb D_4$. Clearly, for $Q$ of type
$\Bbb D_4$, we have $e_Q(3) = 3$. Enlarging one of the branches by an additional vertex, say $\omega$,
the enlarged quiver has precisely one connected subquiver with 3 vertices such that $\omega$ is one
of these vertices. 
	
Such a quiver $Q$ with precisely one branching vertex and three arms is of the form $\Bbb T_{pqr}$ with
$p\ge q\ge r\ge 2$ (the graph $\Bbb T_{pqr}$ is obtained from the disjoint union of graphs of type
$\Bbb A_p, \Bbb A_q, \Bbb A_r$ by identifying one of the end points of these three graphs).
If $Q$ is in addition (as we assume) representation-infinite, then 
$Q$ has a subquiver $Q'$ of type $\widetilde{\Bbb E}_m$ with
$m=6,7,8.$ More precisely: 
We have $r\ge 3$ if and only if $\widetilde{\Bbb E}_6$ is a subquiver of $Q$.
We have $r = 2$ and $q \ge 4$ if and only if $\widetilde{\Bbb E}_6$ is not a subquiver of $Q$,
but $\widetilde{\Bbb E}_7$ is a subquiver of $Q$, and finally $r=2, q=3$ and $p\ge 6$
if and only if neither $\widetilde{\Bbb E}_6$ nor $\widetilde{\Bbb E}_7$ is a subquiver of $Q$,
but $\widetilde{\Bbb E}_8$ is a subquiver of $Q.$
	
Let us exhibit for the quivers $Q$ of type $\widetilde{\Bbb E}_6,\ \widetilde{\Bbb E}_7,\
\widetilde{\Bbb E}_8$
the numbers $e_{Q}(s)$ with $s$ sufficiently small:
$$
\hbox{\beginpicture
\setcoordinatesystem units <1cm,1cm>
\multiput{} at 0 .5  0 -3.5 /
\put{$Q$} at -.35 -.2
\put{$s$} at .1 0.25 
\put{$\widetilde{\Bbb E}_6$} at -.2 -1
\put{$\widetilde{\Bbb E}_7$} at -.2 -2
\put{$\widetilde{\Bbb E}_8$} at -.2 -3
\put{$1$} at 1 0 
\put{$2$} at 2 0 
\put{$3$} at 3 0 
\put{$4$} at 4 0 
\put{$5$} at 5 0 
\put{$6$} at 6 0 
\put{$7$} at 7 0 
\multiput{$6$} at 2 -1  3 -1  /
\multiput{$7$} at 1 -1  4 -1  2 -2  3 -2  4 -2  /
\multiput{$8$} at 1 -2    2 -3  3 -3  4 -3   5 -2  5 -3   6 -3  /
\multiput{$9$} at 1 -3  7 -3    /
\plot -0.6 -.5  7.5 -0.5 /
\plot 0.4 .5  0.4 -3.5 /
\plot 0.4 -0.5  -.5 .5 /
\setshadegrid span <.5mm>
\vshade 3.5  -3.5 -.5 <,z,,> 4.5  -3.5 -.5 <z,z,,> 4.55  -3.5 -1.5 <z,z,,>  5.5  -3.5 -1.5 
 <z,z,,> 5.55  -3.5 -2.5 <z,,,>  7.5  -3.5 -2.5   /
\endpicture}
$$
The numbers in the shaded area still have to be verified, but this is an easy exercise left
to the reader. 

If $Q$ is of type $\Bbb T_{pqr}$ with $r\ge 3$, then $e_Q(4) = n$ (and therefore
$t_Q = 4$). We show this by induction, 
starting with $Q$ of type $\Bbb T_{333} = \widetilde{\Bbb E}_6.$ 
Assume that $Q' \subset Q$ both are quivers of type $\Bbb T_{pqr}$ with $r\ge 3$ such that $Q'$ is
obtained from $Q$ by deleting one vertex, say $\omega$ (and the arrow which involves $\omega$),
then there is just one exceptional representation $M$ of $Q$ of length $4$ with $\omega$ in the
support of $M$ (the support of $M$ is a subquiver of type $\Bbb A_4$). It follows that 
$e_Q(4) = e_{Q'}(4)+1.$ Thus, the induction shows that $e_Q(4) = n$.

If $Q$ is of type $\Bbb T_{pq2}$ with $q\ge 4$, then we claim that
$e_Q(4) = n-1$ and $e_Q(5) = n$ (and therefore  $t_Q = 5$). Again we use induction, 
now starting with the quivers of type $\Bbb T_{442} = \widetilde{\Bbb E}_7.$ 
Assume that $Q' \subset Q$ are quivers of type $\Bbb T_{pq2}$ with $q\ge 4$ such that $Q'$ is
obtained from $Q$ by deleting again one vertex $\omega$.
There is just one exceptional representation $M$ of $Q$ of length $4$ and also just
one exceptional representation $M'$ of $Q$ of length $5$
such that $\omega$ is in the support of $M$ and of $M'$ 
(the support of $M$ and $M'$ are subquivers of type $\Bbb A_4$ and $\Bbb A_5$,
respectively). It follows that $e_Q(s) = e_{Q'}(s)+1,$ for $s=4$ and $s=5$. Thus, the induction
shows that $e_Q(4) = n-1$ and $e_Q(5) = n$.

Finally, let $Q$ be of type $\Bbb T_{p32}$ with $p\ge 6$, then we claim that
$e_Q(s) = n-1$ for $4\le s \le 6$ and $e_Q(7) = n$ (so that $t_Q = 7$).
Also here we use induction, starting now with $Q$ of type $\Bbb T_{632} = \widetilde{\Bbb E}_8.$ 
Assume that $Q' \subset Q$ are quivers of type $\Bbb T_{p32}$ with $p\ge 6$ such that $Q'$ is
obtained from $Q$ by deleting the vertex $\omega$.
Then, for $4\le s \le 7$ there is just one exceptional representation $M$ of $Q$ of length $s$ 
such that $\omega$ is in its support (the support of $M$ is a subquiver of type $\Bbb A_s$), 
thus $e_Q(s) = e_{Q'}(s)+1,$ for $4\le s \le 7$. The induction shows that
$e_Q(s) = n-1$ for $4\le s \le 6$ and $e_Q(7) = n$. 
	
It remains to be seen that 
all indecomposable representations of length at most $t_Q$ are exceptional. Thus, 
let $M$ be an indecomposable representation of $Q$ of length $s \le t_Q$, in particular,
we have $s\le 7$. Assume that $M$ is not exceptional. 
Since $Q$ is a tree quiver, we must have $s \ge 6.$ Also, 
if $s = 6$, then $Q$ has a subquiver of type 
$\widetilde{\Bbb D}_4$, whereas for $s = 7$, the quiver $Q$ has a subquiver of type 
$\widetilde{\Bbb D}_5$. But if $Q$ has a subquiver of type $\widetilde{\Bbb D}_m$ for some $m$,
then $t_Q = 3$ and therefore $s > t_Q$, a contradiction. 
This completes the proof. 
	\medskip
{\bf Remark.} Let $Q$ be a quiver of type $\Bbb T_{pqr}$ with $p\ge q\ge r \ge 2.$ 
Then the proof shows that for $Q$ being representation-infinite, 
most of the indecomposable representations of length at most $t_Q$ are thin.  
For any quiver $Q$ of type $\Bbb T_{pqr}$, there is a unique indecomposable representation of
length $5$ which is not thin. For $Q$ of type $T_{pq2}$ with $q\ge 3$, there are precisely
$2$ indecomposable representation of length $6$ which are not thin. Finally, for $Q$ of type
$\Bbb T_{p32}$ with $p\ge 4$, there are precisely $4$ indecomposable representation of length $7$ which are not thin. 
Thus, if $Q$ is a quiver of type $\Bbb T_{pqr}$, let us denote by $e'_{pqr}(t)$ 
the number of thin exceptional representations of $Q$ of length $t$ and by
$e''_{pqr}(t)$ the number of exceptional representations of $Q$ of length $t$ 
which are not thin.  The relevant numbers for our considerations are the following: 
\vglue-1truecm
$$
\hbox{\beginpicture
\setcoordinatesystem units <.8cm,.7cm>
\plot -1 -1.5  15 -1.5 /
\plot 3 -.5  3 -3.5 /
\plot 7 -.5  7 -3.5 /
\plot 11 -.5  11 -3.5 /
\put{$e'_{pqr}(5)$} at 4 -0.9
\put{$e''_{pqr}(5)$} at 6 -0.9
\put{$e'_{pqr}(6)$} at 8 -0.9
\put{$e''_{pqr}(6)$} at 10 -0.9
\put{$e'_{pqr}(7)$} at 12 -0.9
\put{$e''_{pqr}(7)$} at 14 -0.9

\put{$q\ge 4,\ r=2$} at .7 -2
\put{$p\ge 6,\ q= 3,\ r=2$} at .7 -3

\put{$n-1$} at 4 -2
\put{$n-2$} at 4 -3
\put{$1$} at 6 -2
\put{$1$} at 6 -3

\put{$n-3$} at 8 -3
\put{$2$} at 10 -3

\put{$n-4$} at 12 -3
\put{$4$} at 14 -3
\endpicture}
$$

	\bigskip 
{\bf Question.} It would be of interest to obtain a proof of the implication
(ii) $\implies$ (i) of the following kind: to use the $n$ isomorphism classes of indecomposable representations with fixed length 
in order to construct say arbitrarily large indecomposable representations, or a one-parameter
family of indecomposable representations.
	\bigskip 
{\bf References.}

\frenchspacing
\item{[G]} Gabriel, P.: Unzerlegbare Darstellungen I. manuscripta mathematica 6 (1972), 71-103.
\item{[H]} Humphreys, J\. E\.: Reflection Groups and Coxeter Groups.
   Cambridge Studies in Advanced Mathematics 29. Cambridge University Press (1990).
\item{[K]} Kostant, B.: The principal three-dimensional subgroup and the Betti numbers
   of a complex simple Lie group. Amer.~J.~Math. 81 (1959), 973-1032.
	    
\bigskip\medskip 
{\rmk 
\noindent
Shanghai Jiao Tong University, Shanghai 200240, P. R. China, and \par 
\noindent
King Abdulaziz University, PO Box 80200,  Jeddah, Saudi Arabia. \par
\noindent
e-mail: \ttk ringel\@math.uni-bielefeld.de \par}

\bye